	\documentclass[12pt]{amsart}
\usepackage{amsfonts,amssymb,latexsym,epsfig,amsmath}
\usepackage{amssymb}
\usepackage{amsmath}
\numberwithin{equation}{section}

\newtheorem{prop}{Proposition}[section]
\newtheorem{theor}{Theorem}[section]
\newtheorem{lemma}{Lemma}[section]

\def\qed{~\hfill$\square$}
\def\ep{\varepsilon}

\newcommand{\n}{\hbox{$\scriptstyle\hbox{\rm l\kern-0.22em N}$}}

\renewcommand{\r}{{\mathbb R}}
\def\real{{\mathbb R}}
\newcommand{\R}{\r^{N}}

\renewcommand{\O}{\mathcal{O}}
\newcommand{\e}{\varepsilon}
\newcommand{\tr}{\hbox{\rm tr}}

\def\PS #1 #2{\langle #1, #2 \rangle}

\def\1{1\!\!1}

\def\Ob{{\overline \O}}
\def\dO{{\partial \O}}
\def\sym{S^N}

\def\Z{\mathbb Z}
\def\Kb{{\overline K}}
\def\ub{\bar u}

\def\ut{\tilde u}
\def\zt{\tilde z}

\def\wb{\bar w}

\def\e{\varepsilon}
\def\control{(\alpha_t)_t}
\def\Ex{I\!\!E_x}
\def\Oe{\O_\e}
\def\dOe{\dO_\e}

\def\Fb{\bar F}
\def\ue{u^\varepsilon}

\def\xep{x_\varepsilon}
\def\xoe{\varepsilon^{-1}x}
\def\uz{\tilde z}

 1

\begin{document}

\title[]{\bf 
Ergodic problems and periodic homogenization 
for fully 
nonlinear equations in half-space type domains with Neumann boundary 
conditions}

\author[]
{G. Barles$^{(1)}$, 
F. Da Lio$^{(2)}$, P.-L. Lions$^{(3)}$
and 
P.E. Souganidis$^{(4)}$}

\addtocounter{footnote}{1} \footnotetext{Laboratoire de
Math\'ematiques et Physique Th\'eorique (UMR 6083).  Universit\'e de Tours.
Facult\'e des Sciences et Techniques, Parc de Grandmont, 37200
Tours, France. E-mail address: barles@lmpt.univ-tours.fr }

\addtocounter{footnote}{1} \footnotetext{Dipartimento di Matematica Pura e Applicata,   Universit\`a di
Padova, Via Belzoni 7, 35131 Padova. Italy. E-mail address: dalio@math.unipd.it}

\addtocounter{footnote}{1} \footnotetext{Coll\`ege de France, 11 Place Marcelin Berthelot, Paris 75005 and
CEREMADE, Universit\'e IX Paris-Dauphine, Place du Mar\'echal de Lattre de Tassigny, Paris
Cedex 16, France. E-mail address: lions@ceremade.dauphine.fr}

\addtocounter{footnote}{1} \footnotetext{Department of Mathematics, The University of Texas at Austin,
1 University Station C1200, Austin, TX 78712-0257, USA.
E-mail address: souganid@math.utexas.edu}

\date{}

 \begin{abstract}
We study periodic homogenization problems for second-order pde
in half-space type domains with Neumann boundary conditions. 
In particular, we are interested in ``singular problems'' for which 
it is necessary to determine both 
the homogenized equation and boundary conditions. 
We provide new results for fully nonlinear equations and boundary conditions. 
Our results
extend previous work of Tanaka in the linear, periodic setting 
in half-spaces parallel to the axes of the periodicity, 
and of Arisawa in a rather restrictive nonlinear periodic framework. 
The key step in our analysis is the study of 
associated ergodic problems in domains with similar structure. 
\end{abstract}
 
\maketitle
 
\section{Introduction}

We study issues related to the homogenization and ergodic problems for 
fully nonlinear, non-divergence form, elliptic and parabolic boundary value 
problems in half-space type domains with possibly nonlinear 
Neumann boundary conditions.
In particular, we are interested in problems 
for which it is necessary to identify both the homogenized equation 
and boundary conditions. 
Our results represent a first step towards the resolution of 
the problem in general domains,   which remains open even for 
the linear non-divergence form problem. 
The situation is, of course, different for 
divergence form equations, since, in that case
the boundary conditions are encoded in the variational formulation. 

In order to be more specific and to describe the problem in a clear way, 
we first discuss heuristically two model cases involving linear elliptic 
equations which lead to different behaviors and difficulties. 
To simplify the presentation, we assume that all the functions appearing  
in the examples below have  the needed regularity properties and are 
periodic with respect to the  fast variable $x/\e$, which we 
denote throughout the paper by $y$.

The first problem, which we call  ``regular'', is 
\begin{equation}\label{eq1.1}
\left\{
\begin{array}{rcl}
- \tr (A(x,\xoe)D^2 \ue) - b(x,\xoe)\cdot D \ue+ \ue & = & f(x, \xoe)  \quad \hbox{in  }\O\, ,\\
\displaystyle Du^\e -\gamma & = &g 
\quad \hbox{on  }\dO \; ,
\end{array}
\right.
\end{equation}
where $\O$ is a domain in $\R$, not necessarily of half-space type, 
$D\ue$ and $D^2\ue$ denote respectively the 
gradient and Hessian matrix of the solution $u^\ep$, and
$\gamma$ and $g$ depend only on $x$. 

The typical ``singular'' problem is 
\begin{equation}\label{eq1.2} 
 \left\{
 \begin{array}{rcl}
 - \tr (A(x,\xoe)D^2 \ue) - \displaystyle {\e}^{-1}b(x,\xoe)\cdot D \ue +  \ue 
 & = & f(x, \xoe)   \quad \hbox{in  }\O\, ,\\
\displaystyle Du^\e \cdot \gamma & = & g\quad \hbox{on  }\dO\; ,
 \end{array}
 \right.
\end{equation} 
 which has as 
a particular case  the divergence form equation
$$ - {\rm div}( A(x,\xoe) D\ue) + \ue = f(x,\xoe) \quad \hbox{in  } \O\; .$$

If we use the formal expansion
$$ \ue (x) = \ub (x) + \e v (x, \xoe) + \e ^2 w (x,\xoe) + O(\e^3)\ ,$$
then, in the ``regular'' case, we find that the leading 
(the $\e^{-1}$) term yields 
 $$ - \tr( A(x,y) D_{yy}^2 v ) = 0 \quad\hbox{in  } \R\; ;$$
a standard Liouville-type property implies that $v$ 
is independent of $y$. 
Moreover, when $\gamma$ and $g$ depend only on $x$ and not on $y$, the 
function $\ub$ is expected to satisfy the same boundary condition as 
$\ue$, since the formal expansion yields  
$$ Du^\e \cdot \gamma = D\bar u \cdot\gamma  
+ \e D_x v \cdot \gamma + \e D_y w \cdot \gamma + O (\e^2)\ .
$$

This problem can be treated, at least for uniformly elliptic $A$,
using the perturbed test-function method of Evans \cite{E1}.  
If, however, either $\gamma$ or $g$ depend on $y$, 
different arguments are needed.

The same expansion, in the ``singular'' problem leads ot 
$$\tr ( A (x,y) D^2_{yy} v) + b(x,y) \cdot (D\bar u+ D_y v) = 0\ ,$$
where now the function $v$ actually depends on $y$ in general. 
This interferes with the boundary condition, since the expansion is now
$$ Du^\e \cdot\gamma = D\bar u\cdot \gamma 
 + D_y v \cdot \gamma +  O (\e )\  ,
$$
and clearly 
$v$ plays a role in determining  the boundary condition for 
the limiting equation.

This is the main problem in the singular case.
The key issues are the identification of the  equation 
and the boundary condition for $\bar u$.
In the ``regular'' case, the issue is 
the study of the asymptotic limit for $y$-dependent $\gamma$ and $g$.  

There exists an extensive body of work dealing with 
 the homogenization of the Dirichlet problem for fully nonlinear 
first- and second-order partial differential equations in periodic, 
quasi-periodic, almost periodic and, more recently, stationary ergodic media.
Listing references is beyond the scope of this paper.

Little is, however, known for the Neumann problem except for 
divergence form equations with the usual (co-normal) boundary 
conditions, which are treated 
in the classical book of Bensoussan, Lions and Papanicolaou \cite{BLP}.

In \cite{HT}, Tanaka considered the two model problems discussed
earlier by purely 
probabilistic methods in half-spaces. 
In \cite{A3}, Arisawa studied special cases of homogenization problems again 
in half-space type domains under rather restrictive assumptions.
Our methods are inspired from \cite{A3} but yield more general results.

To study the asymptotic behavior of the $u^{\e}$'s, 
we first consider the usual cell (ergodic) problem which is supposed to give 
the equation inside $\O$.
In the ``singular'' case, it is formulated as follows:

For each $x,p\in\R$, find a unique constant  $\bar\lambda (p,x)$ 
such that there exists a bounded solution $v$ of  the equation
\begin{equation}\label{eq1.3}
- \tr( A(x,y) D_{yy}^2 v ) - b(x,y)\cdot (p + D_y v) = \bar \lambda (p,x) \quad\hbox{in  } \R  \ .
\end{equation}

The map $p \mapsto \bar \lambda (p,x)$ is clearly linear. 
Moreover, 
the only interesting case is when $\bar \lambda (p,x) = 0$, otherwise we 
get a trivial first-order equation for $\ub$. 
We do not discuss here the type of conditions on 
$A$ and $b$ which yield $\bar \lambda (p,x)\equiv 0$. 
We just point out, however, that, 
for divergence form equations, this is true, as 
it can be easily seen by integrating the equation over a period.

For  the boundary behavior, we restrict, for simplicity, to  the half-space 
case. 
To fix the notation, we assume that $\O = \{x_N >0\}$, in which case  
$\O = \e^{-1} \O$.	
The fact that both $\O$ and $\partial\O$ are invariant under the scaling 
is a property which is needed, at least for our approach, in the 
study of general domains.

The natural ergodic problem on the boundary is to find, for each $p,x\in \r^N$, 
a unique constant $\bar \mu (x,p)$ such that there exists a bounded solution 
$w$ of  the boundary value problem 
\begin{equation}
\label{model-ep}
\left\{
\begin{array}{cl}
- \tr( A(x,y) D_{yy}^2 w ) - b(x,y)\cdot (p + D_y w) 
= 0 & \hbox{in  }\O\, ,\\
\noalign{\vskip6pt}
(p + D_y w)\cdot \gamma (x,y) = g(x,y) + \bar \mu (x,p) & \hbox{on  }\dO\  .
\end{array}
\right.
\end{equation}

This type of problem was first studied in \cite{A3} both in 
bounded and in half-space type domains for 
Hamilton-Jacobi-Bellman type equations with oblique Neumann boundary 
conditions but under rather restrictive assumptions. 
We provide here existence and uniqueness results for 
general nonlinear equations and some boundary conditions.

The paper is organized as follows: 
In Section~2 we introduce the hypotheses and 
state and prove the existence  result for the nonlinear 
boundary ergodic problem.  
Section~3	
is devoted to the study of 
the uniqueness properties of $\bar \mu (x,p)$. 
In Section~4,	
we consider homogenization problems.
We conclude with some remarks about the non-periodic case.

\section{Existence results for ergodic boundary problems in half-spaces}
\label{sec:exist}
 
We study  the following boundary ergodic problem~: 
Find a constant $\mu$ such that there exists a bounded, continuous solution of
\begin{equation}\label{eq2.1}
\begin{cases}
F(D^2u,Du,x)&= \lambda\  \mbox{ in $\O$,} \\
L(Du,x) &= \mu\  \mbox{ on $\dO$.} 
\end{cases}
\end{equation}
Here $\O\subset\R$ is a smooth, periodic, half-space-type domain 
(see (O1) below for the definition),
$F$ and $L$ are real valued, 
at least, continuous functions defined on 
$\sym \times \R \times \O$ and $\R\times \dO $ respectively, 
where $S^N$ is the space of $N\times N$ symmetric matrices,  
and $\lambda$ is a constant which is chosen in a suitable way below.  

We formulate next the assumptions on $F$, $L$ and $\O$. 
Note that, for the sake of the clarity of the exposition, in this note,  
we do not attempt to state the most general conditions.

As far as $F$ is concerned, we assume that:  

\medskip
\centerline{{(F0)}\hfill
$F : \sym  \times\R \times\R \to\r$ 
is a continuous, $\Z^N$-periodic function in $x$,\hskip.62truein}
\medskip

\leftline{(F1)\qquad
$\left\{\begin{array}{l}
\text{$F$ is locally Lipschitz continuous 
and there exists a constant $K>0$}\\
\noalign{\vskip6pt}
\text{such that, for all $x,y\in \R$, $p,q \in \R$, and $M,N \in \sym$,}\\
\noalign{\vskip6pt}
 |F(M,p,x) - F(N,q,y)|\leq K
\{ |x-y|(1+|p|+|q| + ||M||+||N||) \\
\noalign{\vskip6pt}
\hskip2.2truein +|p-q| + ||M - N||\}\ ,
\end{array}\right.$ }

\medskip
\leftline{(F2)\qquad
$\left\{\begin{array}{l}
\text{there exists $\kappa >0$ such that, for all $x\in \R$, 
$p\in \R \setminus \{0\}$,}\\
\noalign{\vskip6pt}
\text{and $M,N \in \sym $ with $N\geq 0$,}\\
\noalign{\vskip6pt}
F(M+N,p,x) - F(M,p,x) \leq -\kappa(N\hat p,\hat p) + o(1)||N|| \; ,\\
\noalign{\vskip6pt}
\text{with $o(1)\to 0$ as $|p|\to\infty$,}\\
\end{array}\right.$}
\medskip
\noindent
where, for  $p \in\real^N\setminus \{0\}$, 
$\hat p= p/|p|$.  
\medskip

As far as the domain is concerned, we assume that 
\vskip9pt

\leftline{$(O1)$\hskip1.5 truein 
$\O$ is a smooth, periodic half-space type domain,}
\vskip9pt

\noindent i.e., that:		

\noindent (i) $\dO$ is a $C^2$--boundary, 

\noindent (ii) there exists $(N-1)$ linearly independent vectors 
$f_1, \cdots f_{N-1} \in \Z^N$, 
such that, for all $x\in \dO$,  $\bar x\in \O$ and 
$z_1, \cdots z_{N-1} \in \Z$,
$$ x + z_1 f_1 + \cdots + z_{N-1}f_{N-1} \in \dO\quad \text{ and } 
\quad\bar x + z_1 f_1 + \cdots + z_{N-1}f_{N-1} \in \O\; ,$$
and 

\noindent (iii) there exists a unit vector $f_N \in \R$, orthogonal to 
$f_1, \cdots, f_{N-1}$, and $\bar R >0$ such that, for all 
$x \in \O$ and $R \geq \bar R$, 
$$ x + R f_N \in \O\ .$$
\medskip

Examples of domains satisfying $(O1)$ are 
$$\O = \{ x=(x_1, \cdots, x_N): \; x_N > 0 \}$$ 
or, more generally,   
$$\O =\{ x=(x_1, \cdots, x_N): \; x_N > \psi (x_1, \cdots, x_{N-1}) \},$$ 
where $\psi\in C^2 ({\mathbb R}^{N-1})$ is $\Z^{N-1}$--periodic .
\medskip
 
We denote by $d$ the sign-distance function to $\dO$, normalized to 
be positive in $\O$ and negative in $\R \setminus \Ob$, 
and we recall that, for all $x\in \partial \O$, 
the outward normal $n(x)$ to $\partial \O$ at $x$ is given by 
$$-n(x) = Dd (x)  \ .$$ 

A key ingredient in the existence proof are the up to the boundary 
$C^{0,\alpha}$--regularity results for Neumann boundary value problems
which were obtained by Barles and Da Lio in \cite{BDL}. 
To use them, it is necessary to introduce the following 
set of assumptions on $L$. 

\medskip
\leftline{{(L1)}\qquad 
$\left\{\begin{array}{l}
\text{There exists $\nu>0$ such that, for all $(x,p)\in
\dO \times\R$ and $t >0$,}\\ 
\noalign{\vskip6pt}
\hskip1.3truein 
L(p+t n(x),x)-L(p,x)\ge \nu t\,.
\end{array}\right.$}
\bigskip

\leftline{{(L2)}\qquad 
$\left\{\begin{array}{l}
\text{There exists constant $\Kb >0$ such that, 
for all $x,y\in \dO$ and $p,q \in \R$,}\\ 
\noalign{\vskip6pt}
\hskip.8truein 
|L(p,x)-L(q,y)|\le \Kb \left[(1+ |p|+|q|)|x-y|+|p-q| \right]\,.
\end{array}\right.$}
\bigskip

\leftline{{(L3)}\qquad 
$\left\{\begin{array}{l}
\text{There exists a  $L_\infty \in C(\real^N\times \Ob)$ such that,}\\
\noalign{\vskip6pt}
\text{as  $t\to + \infty$ and locally uniformly in $(p,x)$,}\\
\noalign{\vskip6pt}
\hskip1.7truein 
t^{-1} L(tp,x)\to L_\infty (p,x)\ .
\end{array}\right.$}
\medskip

The regularity results of \cite{BDL} depend on whether $L$ is 
linear or nonlinear and require some additional assumptions on $F$.
In the nonlinear case, 
it is necessary to strengthen (F2) and  to ask for $F$ to be a uniformly 
elliptic, i.e., to  assume that 

\medskip
\leftline{{(F3)}\enspace 
$\left\{\begin{array}{l}
\text{there exists $\kappa >0$ such that, for all $x\in \Ob$, $p\in \R$, }
\text{and $M,N \in \sym $ with $N\geq 0$,}\\
\noalign{\vskip6pt}
 F(M+N,p,x) - F(M,p,x) \leq -\kappa{\tr}(N)\; .\\
\end{array}\right.$}
\medskip

Moreover, we need to assume that    

\medskip
\leftline{(F4)\qquad
$\left\{\begin{array}{l}
\text{there exists $F_\infty \in C(\sym \times \real^N\times \Ob)$ 
such that, locally uniformly in $(M,p,x)$, }\\
\noalign{\vskip6pt}
\text{as $t\to\infty$,\quad
$t^{-1} F(tM, tp,x)\to F_\infty (M,p,x) $}
\end{array}\right.$}

\medskip
For linear boundary condition, i.e., if 
$$L(p,x) = p\cdot \gamma (x) - g(x)\ ,$$
with $\gamma\in C^{0,1}(\partial\O;\real^N)$ 
and $g\in C^{0, \beta}(\partial\O)$ for 
some $\beta\in (0,1)$, instead of (F3) we assume that

\medskip
\leftline{(O2)\qquad
$\left\{\begin{array}{l}
\text{there exists 	
$A\in C^{0,1}(\Ob;S(N))$ and $c_0>0$, such that}\\
\noalign{\vskip6pt}
 A\ge c_0Id\ ,\quad \hbox{and}\quad   A(x) \gamma (x) =n (x)\ \text{ on }\ \dO\ , 
 \end{array}\right.$}
 \medskip
\noindent
and, for all $R>0$, 
\medskip

\leftline{(F5)\qquad 
$\left\{\begin{array}{l}
\text{there exist $L_R, \lambda_{R}> 0$ 
such that, for all $x\in \Ob,$ $|u|\le R,$
$|p|>L_{R}$}\\
\noalign{\vskip6pt}
\text{and $M,\tilde M\in \sym$ with $\tilde M \geq 0$,}\\
\noalign{\vskip6pt}
\text{$F(M+\tilde M,p,x )-F(M,p,x) \le -\lambda_{R} (\tilde M
\widehat{A^{-1}(x)p}, \widehat{A^{-1}(x)p}) + o(1)\|\tilde M\| \, ,$}
\end{array}\right.$}
\medskip
\noindent
where $o(1)\to 0$ as $|p| \to \infty$.

We call (F5) ``adapted  ellipticity'',
since the ellipticity condition is adapted in 
the direction of the oblique derivative through $A$.

We remark that all the assumptions on $F$ are satisfied by 
Hamilton-Jacobi-Bellman or Isaacs Equations under the usual conditions on the 
coefficients, while  (L1),  (L2), (L3) 
are natural conditions for (nonlinear) Neumann boundary conditions.

The constant $\lambda$ in \eqref{eq2.1}  is given by the next proposition.
Since the proof follows the argument of Barles and Souganidis \cite{BS2}, 
we omit it. 

\begin{prop}\label{eqn:whole}
Assume {\rm (F0), (F1)} and {\rm (F2)}. 
There exists a unique constant $\lambda$ such that there exists
a  periodic solution $\ub \in C^{0,1} (\R)$ of
\begin{equation}\label{eq2.2}
F(D^2\ub, D \ub, x) = \lambda \quad \hbox{in  } \R \ .
\end{equation}
\end{prop}

Our result is 

\begin{theor}\label{erghs} 
Assume  {\rm (O1),(F0)} and either {\em (F1), (F3), (F4)}
and {\em (L1), (L2), (L3)} if  $L$ is nonlinear or 
{\em (F1), (F2), (F5), (O2)}, $\gamma\in C^{0,1}(\partial \O)$
and $g\in C^{0,\beta} (\partial \O)$, if $L$ is linear. 
There exists $\mu$ such that \eqref{eq2.1}
has a continuous bounded viscosity solution, which has the same 
periodicity property as the domain $\O$.
\end{theor}

\medskip
\noindent {\bf Proof.} 
We concentrate first on the nonlinear case. 
To simplify the exposition, we assume that $f_N = e_N$ 
and $\text{span}(f_1, \cdots, f_{N-1}) = \text{span}(e_1, \cdots, e_{N-1})$, 
where $(e_1, \cdots , e_N)$ is the standard orthonormal basis of $\R$. 
In addition, we assume that $0 \in \dO$.

For $0<\varepsilon < \alpha < 1$, we introduce the approximate problem
\begin{equation}\label{eq2.3}
\left\{ \begin{array}{r}
F(D^2\ut,D\ut,x)+\e \ut=\lambda + \e \ub  \mbox{ in $\O$,} \\ 
L( D\ut,x) +\alpha\ut=0 \mbox{ on $\dO $,} 
\end{array}\right.
\end{equation}
where $F$ and $\bar u$ are  given by Proposition~2.1.

The existence and uniqueness of $\ut$ follows from classical arguments 
from the theory of viscosity solutions, and in particular, 
the Perron's method \cite{I2} and the comparison results \cite{b2}. 
Moreover, in view of its uniqueness, 
the solution has the same periodicity properties as the domain, 
i.e., for  all $x\in \Ob$ and $k_1 \cdots k_{N-1} \in \Z$,  
$$ \ut (x + k_1 e_1 + \cdots + k_{N-1}e_{N-1}) = \ut (x)\, .$$

Finally standard comparison arguments (see \cite{b2})  yield that 
$$ \max_{\Ob}(\ut - \ub) \leq \alpha^{-1} 
\sup_{x\in \R, |e|\leq ||D\ub||_\infty}\, |L(e,x)|\; ,$$
which, of course, implies that 
$$\alpha \ut \text{ is uniformly bounded in } x,\alpha \text{ and }\ep\ ;$$
notice that the $\varepsilon$-term in \eqref{eq2.3}  is  
introduced just in order to prove this estimate.

For $R >0$ sufficiently large, we consider next the domain 
$$\O_R = \{x \in \O: \; x_N < R\}$$
and  solve the new problem
\begin{equation}\label{eq2.4}
\left\{ \begin{array}{rl}
F(D^2\ut_R,D\ut_R,x)+\e \ut_R=\lambda + \e \ub & \mbox{in $\O_R$,}\\ 
\noalign{\vskip6pt}
D\tilde u_R\cdot n_R = 0 & \mbox{on $\{x_N = R\}$,} \\
\noalign{\vskip6pt}
L(D\ut_R,x) +\alpha\ut_R=0& \mbox{on $\dO $,} 
\end{array}\right.
\end{equation}
where $n_R$ is the external normal to $\O_R$ on $\{x_N = R\}$. 

Again it is a standard fact that \eqref{eq2.4} admits a unique continuous 
viscosity solution $\tilde u_R$ 
with the same periodicity properties as the domain. 
%
Moreover, the maximum principle  yields that, for all $x\in \Ob$,
$$ |\ut_R (x) - \ut_R(0)| \leq || \ut_R - \ut_R(0)||_{L^\infty(\dO)} \; ,$$
and 
since $\dO$ is periodic, we may assume that all the $u_R-\tilde u_R(0)$'s  
attain their maximum and minimum values at points  
which remain in a compact subset of $\R$.

Letting $R\to \infty$ and using the standard fact that, as $R\to\infty$,
$\tilde u_R\to \tilde u$ in $C(\O)$, 
we find that, 
for all $x\in \Ob$,
$$ |\ut (x) - \ut (0)| \leq || \ut  - \ut (0)||_{L^\infty(\dO)} \; .$$

We claim next that 
$$ || \ut  - \ut (0)||_{L^\infty(\Ob)} = || \ut  - \ut (0)||_{L^\infty(\dO)}$$
remains bounded as first $\varepsilon\to0$ and then $\alpha  \to 0$.
To this end, we argue by contradiction and introduce the function 
$\tilde w: \overline{\O} \to {\mathbb R}$ given by
$$\tilde w= \frac{\ut - \ut (0)}{|| \ut  - \ut (0)||_{L^\infty(\Ob)}}\; .$$

It is immediate that 
$|\tilde w| \leq 1$, $\tilde w(0)=0$ and $w$ is a solution of a 
nonlinear Neumann type 
problem with nonlinearities satisfying assumptions {(F1)}, {(F4)}, 
{(F5)}, {(L1)}, {(L2)} and {(L3)} with uniform constants.
The $C^{0,\alpha}$--estimates of \cite{BDL} yield that the 
$\tilde w$'s are locally uniformly bounded in  $C^{0,\alpha}$. 
Extracting a subsequence, we may assume  that,  as 
$\varepsilon\to0$ and $\alpha\to0$, the $\tilde w$'s converge to $\bar w$, 
which solves
\begin{equation}\label{eq2.5}
\left\{ \begin{array}{r}
F_\infty (D^2 \wb ,D \wb,x )=0  \mbox{ in $\O$,}\\  
\noalign{\vskip6pt}
L_\infty (D \wb,x )=0  \mbox{ on $\dO$,} 
\end{array}\right.
\end{equation}
and, moreover, 
$$\wb (0) = 0,\ \text{ and }\ ||\wb||_\infty = 1$$
with the sup--norm 
achieved at some point of $\dO$, a fact which 
contradicts the strong maximum principle of 
\cite{BDL-eb}.

Since, by the previous step, $\ut  - \ut (0)$ remains 
bounded, we can use again the $C^{0,\alpha}$-local estimates of \cite{BDL} 
and pass to the limit, letting first $\e \to 0$  and then $\alpha\to0$. 
It follows that, up to subsequences, 	
$$- \alpha  \ut (0) \to \mu\ .$$  		

The proof of the linear case follows  the same arguments using 
a different $C^{0,\alpha}$-local estimate in \cite{BDL}. 
It is worth pointing out that the strong maximum principle of 
\cite{BDL-eb} is still valid under assumption {(F5)}.\qed

\section{Uniqueness results for boundary ergodic problems 
for Hamilton-Jacobi-Bellman type equations in half-space type domains} 
\label{sec:uni}

We present here two proofs for the uniqueness of the boundary 
ergodic cost $\mu$. 
The first one assumes that $F$ is convex and 
uses the stochastic control interpretation of the problem. 
%
It turns out, however, that it is also possible to present a proof which is 
based entirely on  pde-type arguments and does not utilize convexity 
and stochastic control.

\medskip
\noindent {\it 3.1. The stochastic control proof.}

We consider the solution 
$((X_t)_{t\geq 0} ,(k_t)_{t\geq 0})$ 
of the stochastic differential equations
\begin{equation}\label{ROe}	
\left\{ \begin{array}{ll}
dX_t = b(X_t, \alpha_t)dt +\sqrt2 \sigma (X_t,\alpha_t)dW_t
- dk_t\; , & X_0=x \in \Ob , \\  
\noalign{\vskip6pt}
k_t =\int_0^t 1_{\dO} (X_s)
\gamma (X_s)d|k|_s \; ,& X_t \in \Ob \; ,\quad  t\geq 0 \; , 
\end{array} \right.
\end{equation}
where $(X_t)_{t\geq 0}$ is a continuous process in $\R$
and $(k_t)_{t\geq 0}$ is a process with bounded variation. 
Here $(W_t)_{t\geqq 0}$ is an $N$-dimensional Brownian motion, the control 
$(\alpha_t)_{t\geqq 0}$ 
is a progressively measurable process with 
respect to the filtration associated to the Brownian motion with values in  
a compact metric space $\mathcal{A}$  
and the drift $b$ and the diffusion matrix $\sigma$ satisfy 
the classical assumptions

\medskip
\leftline{(C1)\qquad 
$\left\{\begin{array}{l}
\text{$b\in C(\Ob\times \mathcal{A};\mathbb{R}^N)$, 
$\sigma\in C(\Ob\times \mathcal{A};S(N))$  
and, for each $\alpha\in\mathcal{A}$,}\\ 
\noalign{\vskip6pt}
\text{$b(\cdot,\alpha)\in C^{0,1}(\Ob;\mathbb{R}^N)$, 
$\sigma(\cdot,\alpha)\in C^{0,1}(\Ob\times S(N))$ and $a=\sigma\sigma^T$}\\
\noalign{\vskip6pt}
\text{is uniformly elliptic, with all the constants uniform in $\alpha$.}
\end{array}\right.$}
\medskip

It is known that (C1) yields 
the existence of a unique solution 
$((X_t)_{t\geq 0} ,(k_t)_{t\geq 0})$  of \eqref{ROe}. 
(See Lions and Sznitman \cite{LiSz} for the existence  and 
Barles and Lions \cite{BaLi}  for the uniqueness.)

The nonlinearity $F$ of the associated Bellman equation
and the boundary condition for the ergodic boundary problem are given 
respectively by
\begin{equation}\label{eq3.2}
 F(M,p,x)=\sup_{\alpha\in \mathcal{A}}\,
\Big\{- {\rm tr}[
a(x, \alpha)M]-b(x,\alpha)\cdot p -f(x,\alpha)\Big\}\ 
\end{equation}
and 
\begin{equation} \label{eq3.3}
L (Du,u,x) = 
Du\cdot  \gamma - g + \mu =0 \quad \mbox{on $\dO$.}
\end{equation}

We assume that 

\medskip
\leftline{(C2)\qquad 
$\left\{\begin{array}{l}
\text{there exist $\beta\in (0,1)$ and $\nu>0$ 
such that, for all $\alpha \in \mathcal{A}$, 
$f(\cdot,\alpha)$,}\\ 
\noalign{\vskip6pt}
\text{$\gamma$ and $g$ satisfy : $f(\cdot,\alpha) \in BUC (\real^N),\  g\in C^{0,\beta}(\partial\O)$,}\\
\noalign{\vskip6pt}
\text{$\  
\gamma\in C^{0,1}(\partial\O;  \real^N)
\text{ and  }\gamma (x) \cdot n (x) \geqq 0$ on $\partial\O.$}
\end{array}\right.$}
\medskip

Finally, it is necessary to make the 
additional ``ergodic''-type assumption that  
\medskip

\leftline{(E1)\qquad 
$\left\{\begin{array}{c}
\text{there exists a bounded, Lipschitz continuous subsolution $\hat w$ of}\\ 
\noalign{\vskip6pt}
\hat F (D^2\hat w, x,-f_N + D\hat w, x) = 0 \quad\hbox{in  }\R\; ,
\end{array}\right.$}
\medskip

\noindent
where $f_N$ appears in $(O1)$ and 
$$\hat F(X,p,x) = \sup_{\alpha\in \mathcal{A}} 
\Big\{ -\frac12 \text{tr} \left[ a(x,\alpha)X\right] 
- b(x,\alpha)\cdot p\Big\}.$$

We remark that, for 
linear equations,  (E1) is not a real restriction. 
 Indeed, in the Introduction, we explain that 
$\bar\lambda (p,x) = 0$ for all $p$ is natural in 
our framework and the condition (E1) reads $\bar\lambda (-f_N,x) = 0$.

The result is 

\begin{theor}\label{uni}
Assume {\rm (C1), (C2)} and {\rm (E1)}. 
There exists at most one constant 
$\mu $ which  solves  the boundary ergodic control problem 
\eqref{eq2.1} with  $\lambda$ and $F$ and $L$ 
given by Proposition~\ref{eqn:whole} and \eqref{eq3.2} and \eqref{eq3.3} 
respectively. 
\end{theor}

The key step  in the proof of the uniqueness is the following lemma.

\begin{lemma}\label{condonk} 
Assume {\rm (C1), (C2)} and {\rm (E1)}.
There exists $x \in \O$ such that
\begin{equation}\label{condk}
\inf_{\control}\Ex \left[ \int_0^{+\infty} d|k|_s \right] = +\infty\; .
\end{equation}
\end{lemma}

\noindent {\bf Proof.} 
Following the proof of Theorem~\ref{erghs}, we assume that $f_N = e_N$, 
and we work in 
$$\O_R = \{x \in \O\ : \; x_N < R\}\ ,$$

Define 
$$v(x)= - x_N + w(x) $$ 
where $w$ is given by Proposition~2.1. 
Then 
$$ \hat F (D^2 v, Dv, x) = 0 \quad\hbox{in  }\R\ ,$$

Since $v$ is Lipschitz continuous and $w$ is bounded, the exists a 
$C>0$ such that $v$ 
$$ Dv\cdot \gamma \leq C  \quad \mbox{on $\dO$}
\ \text{ and }\ 
 v \leq -R + C  \quad \mbox{on $x_N = R$.}
$$

Let  $\tau_x$ be the first exit time through the boundary  $\{x_N = R\}$ 
of the process $(X_t)_{t\geqq 0}$ starting at $x$. 
Using the dynamic programming principle for the stopping times 
$\tau_x$ we find that, for all $t>0$,
$$ v (x) \leq C \inf_{\control}\,\Ex\left[
\int_0^{t\wedge \tau_x}  d|k|_s \,\right] - (R-C) \chi(x,t)\; ,$$
where 
$$ \chi(x,t)
= \inf_{(\alpha_t)_{t\geqq 0}}\,\Ex\left[\1_{\{\tau_x \leq t\}} \,\right]\ .
$$ 

Since $\sigma$ is nondegenerate, it follows that, for each fixed $R$, 
as $t\to\infty$ and locally uniformly in $x$, 
$$ \chi(\cdot,t) \to 1 \ ,$$
which, in turn, implies that, locally uniformly in $x$, 
$$  \liminf_{t \to \infty}\,  \inf_{(\alpha_t)_{t\geqq 0}}\,\Ex\left[
\int_0^{t\wedge \tau_x}  d|k|_s \,\right] \geq C^{-1}(R -C + v(x))\  .$$

Since $R$ is arbitrary, 
the result now follows for all $x\in \O$.\qed


We return now to the proof of Theorem~3 which is based on the intuitive 
idea that,  to have a unique $\mu$, 
the boundary needs to be seen in a ``sufficient way''. 
This last fact which is quantified in a precise way by  (\ref{condk}).  
A similar fact is also a key point for the uniqueness of 
$\mu$ in bounded domain.
\bigskip

\noindent {\bf Proof of Theorem~\ref{condonk}.}
Assume that there exist two constants $\mu_1$ and $\mu_2$ such that 
$\mu_1 <  \mu_2$ with corresponding solutions $u_1$ and $u_2$. 
The dynamic programming principle
or the uniqueness for the associated time dependent problem  yield, 
for all $t>0$  and $i=1,2$, 
$$
u_1 (x) = \inf_{(\alpha_t)_{t\geqq0}}\,\Ex\left[
\int_0^t [f(X_s,\alpha_s)+\lambda ] ds 
+ \int_0^t [g(X_s)+\mu_1]d|k|_s + u_1 (X_t) \,\right] \ .
$$ 

Choosing an $\e$-optimal control for $u_2$ and subtracting the 
inequalities for $u_1$ and $u_2$ we find 
$$ u_1 (x) -  u_2 (x) \leq (\mu_1 -  \mu_1) 
\Ex\left[ \int_0^t d|k|_s\right] + 
u_1 (X_t) -  u_2 (X_t) + \e \; ,$$
which yields the estimate 
$$(\mu_2 - \mu_1) \Ex\left[ \int_0^t d|k|_s\right] 
\leq 2(||u||_\infty + ||\tilde u||_\infty ) + \e\, . $$

If $\mu_2 >\mu_1$ this last inequality contradicts Lemma~3.1.\qed
\medskip

 

It turns out that, for $\mu$ to be unique, it is necessary to assume (E1).
Indeed, consider the simple Neumann problem 
$$ -\varphi '' - \varphi' = 0\quad \hbox{in  }(0, +\infty)\;  
\ \text{ and }\ 
 - \varphi'(0)  = \mu \ .$$	
Since $\phi (x) = \mu \exp (-x)$ solves this problem for all $\mu \geqq 0$, 
it turns out that since
$\varphi(x)= \mu \exp(-x)$, 
any $\mu\geqq 0$ is a solution of the associated ergodic problem. 

On the other hand, a direct integration yields that 
there does not exist a bounded subsolution of 
$$ -\psi '' - \psi' + 1 = 0\quad \hbox{in  }(0, +\infty)\ . $$

Finally, we remark that in \cite{A3} 
the uniqueness is proved in the uniformly elliptic case under the 
additional assumption that, for all $x$ and $\alpha$,
$$ b_N (x,\alpha) \leq 0 \ ,$$ 
a fact which allows to choose $w\equiv 0$ in (E1).

\medskip
\noindent {\it 3.2. A pde proof and a more general result.}

We state and prove using pde techniques a stronger uniqueness result 
than the one asserted in Theorem~3.1. 
To this end, we need to assume that there exist

\medskip
\leftline{(F6)\qquad 
$\left\{\begin{array}{l}
\text{a homogeneous of degree one, uniformly elliptic,  
$\hat F \in C(\sym \times  \R \times \R)$}\\
\noalign{\vskip6pt} 
\text{such that, for all $ x \in \R$, $p,q  \in \R$ 
and $M_1,M_2 \in \sym$,}\\
\noalign{\vskip6pt}
\text{$\qquad F(M_1,p,x) - F(M_2,q,x) \leq \hat F (M_1-M_2,p-q,x)\; ,$}
\end{array}\right.$}
\medskip

\noindent
and
\medskip

\leftline{(E2)\qquad 
$\left\{ \begin{array}{l}
\text{a Lipschitz continuous subsolution $v$ of}\\
\noalign{\vskip6pt}
{}\hskip1in   \hat F (D^2v, Dv, x) = 0 \quad\hbox{in  }\O\; ,\\
\noalign{\vskip6pt}
\text{such that, uniformly with respect to $x'$,
$v(x',x_N) \to -\infty$ as $x_N \to + \infty$,} 
\end{array}\right. $}
\medskip

\noindent
where, for $x\in\real^N$, we write $x= (x',x_N)$ with $x' = (x_1,\ldots,x_N)$.

It is immediate that the convex function $F$ given by \eqref{eq3.2}
satisfies (F6), while (E2) is a weaker version of (E1).
We also remark that, although all these assumptions appear a bit artificial,
\cite{A3}  provides some counterexamples.

The result is 
\begin{theor}\label{unibis}
Assume {\rm (F6)} and {\rm (E2)}.
If $\hat F$ satisfies {\rm (F1)}, there 
exists at most one constant $\mu$ solving the  boundary 
ergodic control problem.
\end{theor}

\noindent {\bf Proof of Theorem~\ref{unibis}.} 
1. Let $\mu_1$, $\mu_2$, $u_1$, $u_2$ and $\O_R$ be as in the proof of 
Theorem~\ref{uni} and assume that $\mu_2 > \mu_1$. 
Then $w= \tilde u -u$ is a bounded supersolution of 
\begin{equation*}
\left\{
\begin{array}{l}
\hat F (D^2 w, Dw, x) \geq 0 \quad\hbox{in  }\O \; ,\\
Dw\cdot \gamma  \geq  \mu_2 - \mu_1  
\quad \mbox{on $\dO$.}
\end{array}\right.
\end{equation*}

For some $c>0$,  assume that there exists 
$\psi^R \in C(\O\times [0,\infty))$  such that 
\begin{equation}\label{eq:A}
\psi^R \leqq c\text{ on } \{x_N=\real\}\text{ and } 
\psi_R (\cdot,0)\leqq w\text{ on } \partial\O\ ,
\end{equation}
and
\begin{equation}\label{eq:B}
\left\{ \begin{array}{l}
\psi_t^R + \hat F(D^2\psi^R, D\psi^R,x) \leqq 0\text{ on } \O_R\ ,\\
\noalign{\vskip6pt}
\displaystyle\frac{\partial\psi^R}{\partial\gamma} \leqq \tilde \mu-\mu 
\text{ on } \partial\O_R\ .
\end{array}\right.
\end{equation}

The comparison principle of viscosity solutions and the fact that 
$w$ is bounded yield the existence of another constant $C>0$ such that 
$$\psi^R \leqq w +C\ \text{ in }\ \O\times (0,\infty)\ .$$

This last inequality leads to a contradiction, if it is possible to choose 
the $\psi^R$'s so that, in addition to \eqref{eq:A} and \eqref{eq:B},  
they also satisfy,  
\begin{equation}\label{eq:C}
\lim_{R\to \infty} \lim_{t\to\infty} \psi^R (\cdot,t) = +\infty\ 
\text{ locally uniformly in $\O$.} 
\end{equation}

For each $R>0$,  we define  $\psi^R$ by  
\begin{equation*}	
\psi^R  = \delta v + \delta m_R\chi^R \ ,
\end{equation*}
where $\delta >0$ is chosen below sufficiently small enough, 
$$m_R = \inf_{x_N =R} (-v(x',x_N)) \to \infty\ \text{ as }\ R\to\infty\ ,$$
and $\chi^R \in C(\O\times (0,\infty))$ is such 
\begin{equation}	\label{eq:E}
0\leqq \chi^R \leqq 1  
\text{ and } 
\lim_{t\to\infty} \chi^R (\cdot,t) =1\text{ locally uniformly in }\O\ .
\end{equation}

Such $\psi^R$ clearly satisfies \eqref{eq:C} and the first inequality 
in  (11). 
To prove that $\psi^R$ is a subsolution we compute
$$ \psi_t^R + \hat F ( D^2\psi^R , D\psi^R,x) 
= \delta m_R \chi^R_t + \hat F (\delta D^2 v + \delta m_R D^2 \chi^R,
\delta D v + \delta m_R D \chi^R, x)$$
$$ \leq  \hat F (\delta D^2 v, \delta D v,x ) +  \delta m_R 
\left( \chi^R_t - {\mathcal M}_-(D^2 \chi^R) - C |D \chi^R| \right)\; ,$$
where ${\mathcal M}_-$ is the minimal Pucci's operator
associated with $\hat F$.

This computation shows that it is enough to choose 
$\chi^R$ to be the unique solution of 
$$\begin{cases}
\chi^R_t - {\mathcal M}_-(D^2 \chi^R) - C |D \chi_R| = 0 
\ \text{ in }\ \O_R \times (0, +\infty)\ ,\\
\noalign{\vskip6pt}
D\chi^R \cdot \gamma =0\ \text{ on }\ \dO\times (0,+\infty)\ ,\\
\noalign{\vskip6pt}
\chi^R = 1 \ \text{ on }\ \{x_N = R\}\times (0, +\infty),\\
\noalign{\vskip6pt}
\chi^R = 0 \ \text{ on }\ \Ob \times \{0\}\ ,
\end{cases}$$
which satisfies \eqref{eq:E}.

With this choice of $\chi^R$, $\psi^R$ satisfies 
the subsolution inequality. 
The fact that $v$ is bounded yields, if $\delta$ is chosen sufficiently small,
that the boundary condition in \eqref{eq:A} and the second inequality in 
\eqref{condk}  are also satisfied. 

\rightline{\qed}

%
%
%


\medskip
We conclude this section emphasizing that an import issue is to understand
the way on which $\mu$ depends on $F$ and $L$.
This question was addressed in \cite{BDL-eb} in bounded domains.  
In the context of this paper, the dependence on $L$ is an easy consequence 
of the existence proof, while the dependence on $F$ is not so obvious.

\section{The homogenization problem in half-space type domains\\
with oscillating boundaries}\label{sec:hom}

To make the main ideas of our approach clear,
we begin with the linear problem and  
 we consider domains of the form
$$ \Oe = \left\{x\in \R\, :\; x_N > \e \psi (\xoe)\right\} 
= \e^{-1} \O = 
\left\{y\in \R\, :\; y_N > \psi (y ')\right\} \; ,$$
where, as before, for $x\in\real^N$,  
$x' = (x_1, \cdots, x_{N-1})$ and $\psi$ is a smooth (at 
least $C^{2,1}$) $\Z^{N-1}$-periodic function.

The model homogenization problem we are interested in is 
\begin{equation}\label{homlin}
 \left\{
 \begin{array}{rcl}
 - \tr (A(x, \xoe)D^2 \ue) - \displaystyle \e^{-1} b(x,\xoe)\cdot D \ue 
 +  \ue & = & f(x, \xoe)   \quad \hbox{in  }\Oe \, ,\\
 \noalign{\vskip6pt}
Du^\ep \cdot \gamma & = & g(x,\xoe) \quad \hbox{on  }\dOe\; .
  \end{array}
 \right.
\end{equation}

As far as the coefficients are concerned, we assume  that 
$$\text{(H1)}\qquad \left\{ 
\vcenter{\baselineskip=20pt \vsize=50pt\hsize=5.0truein \noindent
$A$, $b$, $f$, $\gamma$ and $g$ are 
bounded, Lipschitz continuous, $\Z^N$-periodic with respect to the 
fast variable and there exists $\nu >0$ such that
$ A \geq \nu Id$.  		
} \right.  \hskip1truein$$

In order to formulate the result, we introduce first the associated  
cell problem which is used for the equation inside the domain. 
We assume that 
$$\text{(H2)}\qquad \left\{ 
\vcenter{\baselineskip=20pt \vsize=50pt\hsize=5.0truein \noindent
for all $x, p\in \real^N$, there exists a bounded, $\Z^N$-periodic in $y$, solution $v(p,x,y)$ of
 \newline
${}\qquad - \tr( A(x,y) D_{yy}^2 v ) + b(x,y)\cdot (p + D_y v) = 0 
\quad\hbox{in  } \R\; ,$
\newline
which depends smoothly on $(p,x,y)$.
} \right.  \hskip.8truein$$
\medskip

The cell problem for the second corrector is 
\begin{equation*}
\text{(H3)}\qquad \left\{ 
\vcenter{\baselineskip=20pt \vsize=50pt\hsize=5.0truein \noindent
for all $x$, $p \in \real^N$ and $M\in S^N$, 
there exists a unique constant $\Fb (M,p, x)$ 
such that the equation \newline
\centerline{$- \tr( A(x,y)\left(M + 2 (D_{xy}^2 v + MD_{yp}^2v)  
+ D_{yy}^2 w \right) - b(x,y)\cdot$}\newline
\centerline{$(D_{x}^2 v + MD_{p}^2v + D_y w) 
= f(x,y) - \Fb (M,p, x) \quad\hbox{in  } \R $} \newline 
has a bounded  solution $w$,			
which depends smoothly on $(M,p,x,y)$.
} \right. \hskip1truein 
\end{equation*}

\medskip
We remark that in this linear context the assumed uniqueness  yields that 
$\Fb$ is an affine function of $p$ and $M$.

As it was mentioned in the Introduction, (H2) and (H3) contain two assumptions. 
One  is the fact that the cell problems have, 
for all $x$, $p$ and $M$,  bounded solutions.
The second is, of course, 
the smoothness of $v$  and $w$ with respect to all the variables. 
This is a technical assumption, which is difficult to avoid 
and, most probably, is true for smoother coefficients. 

We turn next to the boundary condition. 
The results of Sections~\ref{sec:exist} and \ref{sec:uni} 
yield, for each $p,x \in \real^N$, 
the existence of a unique constant $\bar \mu (p,x)$ 
such that the boundary ergodic problem
\begin{equation}\label{erghomlin}
 \left\{
 \begin{array}{rcl}
- \tr( A(x,y) D_{yy}^2 z )-  b(x,y)\cdot (p + D_y z) & = & 0 
\quad\hbox{in  } \O \, ,\\
\noalign{\vskip6pt}
(Dz+p) \cdot \gamma (x,y) 
 & = & g(x,y)-\bar \mu (p,x) \quad \hbox{on  }\dO\; ,
  \end{array}
 \right.
\end{equation}
has a bounded, periodic in $y'$, solution $z$. 

Now we are in position to state our result.

\begin{theor}\label{th:homlin} 
Assume {\rm (H1), (H2)}, {\rm (H3)} and that  
$\bar \mu$ satisfies {\rm (L1)} and {\rm (L2)}. 
The family $(\ue)_{\e>0}$ converges locally uniformly, 
as $\e\to0$, to the unique solution $\bar u$ of
\begin{equation}\label{limeqnL}
 \left\{
 \begin{array}{rcl}
\Fb (D^2 \bar u,D\bar u,x) + \bar u = 0 \quad\hbox{in  } \{x_N >0\} \, ,\\
\noalign{\vskip6pt}
\bar \mu (D\bar u,x) = 0 \quad \hbox{on  } \{x_N = 0\}\; .
  \end{array}
 \right.
\end{equation}
\end{theor}

The uniqueness of $\bar \mu (p,x)$ yields that $\bar \mu$ 
is an affine function of $p$ and, hence,  of the form
$$ \bar \mu (x,p) = \bar \gamma (x)\cdot  p - \bar g (x)\; ,$$ 
for some $\bar\gamma \in C(\partial\O,\real^N)$ and 
$\bar g\in C(\partial \O)$.

The stability and the uniform $C^{0, \alpha}$-estimates also yield that 
the boundary ergodic problem has a solution.
To prove, however,  that $\bar \gamma$ and  $\bar g $ are 
Lipschitz continuous is more 
difficult, since it is necessary to take into account  the dependence an 
$x$ of $A$, $b$ and $f$, i.e., the equation inside the domain. 
This is why we assume that $\bar \mu (p,x)$ satisfies 
(L1) and (L2) in order to have a uniqueness result for (\ref{limeqnL}).

\medskip
\noindent{\bf Proof.} 
Since the proofs that the $u_\varepsilon$'s converge to 
viscosity subsolutions and  supersolutions of (\ref{limeqnL}) are similar, 
here we present the argument only for the subsolution case. 
 
To this end, we introduce the half-relaxed limit 
$$ \ub (x)= {\limsup}^* \ue (x) = \limsup_{x'\to x,\, \e\to 0}
\ue (x')$$
and assume that $\bar x$ is a strict maximum of $\bar u-\phi$, where $\phi$ 
is a smooth test-function.

If $\bar x \in \{x_N >0\}$, the conclusion follows easily  using 
the perturbed test-function.
Indeed consider maximum points $x_\e$ of
$$ \ue (x) - \left( \phi(x) + \e v (D\phi (x), x, \e^{-1}x) 
+ \e^2 w (D^2\phi (\bar x), D\phi (\bar x), \bar x,\e^{-1}x)\right)
\; .$$

The perturbed test-function is smooth and the proof follows as in 
the formal asymptotic expansion. 
It is worth pointing out that technically the smoothness of $v$ is 
required because of its dependence on $D\phi(x)$ and $x$ is important in the 
$\e^{-1}$-term. 
On the contrary, for $w$, the terms are less sensitive and the
dependence on $\bar x, D\phi (\bar  x)$ and $D^2\phi (\bar  x)$ is enough.

This is exactly the difficulty we face, when $\bar x \in \{x_N = 0\}$. 
Assuming that the solution $z$ of (\ref{erghomlin}) is smooth with respect 
to all variables is   too restrictive and, most probably, uncheckable. 
We avoid this in the following way. 
First we write 
$$z (y; D\phi(x),x) = v (x; D\phi (x), y) +\tilde z (y,x)\ .$$

It follows that $\tilde z$ solves the problem  
\begin{equation}\label{erghom2}
 \left\{
 \begin{array}{l}
-\tr( A(x,y) D_{yy}^2 \tilde z) - b(x,y)\cdot D_y \tilde z  
=  0 \quad\hbox{in  } \O \, ,\\
\noalign{\vskip6pt}
D\tilde z\cdot\gamma 
 =  \tilde g(x,y)-\bar \mu (D\phi(x),x) \quad \hbox{on  }\dO\; ,\\
  \end{array}
 \right.
\end{equation}
where 
$$\tilde g(x,y) = g(x,y) - D\phi (x)\cdot \gamma (x,y) 
- D_y v (x, D\phi (x), y) \cdot \gamma (x,y)\ .$$

Since $\tilde z$  is not a smooth function of $x$ and, perhaps, $y$, 
we use Theorem~\ref{erghs} to solve the boundary value problem 
\begin{equation}\label{erghom3}
 \left\{
 \begin{array}{rcl}
 \displaystyle\max_{|x'-x|\le \delta}\left( - \tr( A(x',y) D_{yy}^2 
 \tilde z^\delta) - b(x',y)\cdot D_y \tilde z^\delta \right) 
 & = & 0 \quad\hbox{in  } \O \, ,\\
 \noalign{\vskip6pt}
 \displaystyle \max_{|x'-x|\le \delta}
\left(D\tilde z^\delta \cdot \gamma  -\tilde g(x',y)\right)
& =  & -\bar \mu^\delta \quad \hbox{on  }\dO\; .
  \end{array}
 \right.
\end{equation}

Using the available estimates on $\zt^\delta$ and $\bar \mu^\delta$, it is 
then easy 
to prove that, as $\delta \to 0$, $\zt^\delta$ converges locally uniformly to 
a solution of (\ref{erghom2}) and, more importantly, 
using the uniqueness of $\bar \mu (D\phi(x),x)$,
$\bar \mu^\delta$ that converges to $\bar \mu (D\phi(x),x)$.  

Next we apply the perturbed test-function method in the following way. 
We look at maximum points $x_\varepsilon$ of
$$\ue (x) - \phi(x) - \e v (x, D\phi (x), \xoe) 
-  \zt^\delta (\xoe) 
- \e^2 w(\bar x, D\phi (\bar x), D^2\phi (\bar x), \xoe)
\; .$$ 

To conclude, it suffices to remark that, since $\xep \to \bar x$ as $\e\to0$, 
we have $|x_\e -\bar x|\le \delta$ for $\e$ small enough 
and (\ref{erghom3}) provides the right inequalities to conclude, since the 
maximum is bigger than the value at $x=\xep$. 
The fact that $\uz^\delta$ may be a non-smooth 
viscosity solution of (\ref{erghom3}) creates no real difficulty 
in the argument. 
Indeed, it suffices to double variables using the test-function associated to 
(\ref{erghom3}).~\qed
\medskip

We remark 
that $v$ carries all the information to build the second corrector $w$. 
On the other hand, $\uz$ appears to be necessary 
to treat the boundary condition. 
For this reason, our decomposition of the first corrector appears natural.


We turn next to the nonlinear problem. 
We use the same notations and, in particular, 
the same domain as in the previous section, and, 
to simplify the presentation, we consider a model problem of the form
\begin{equation}\label{homnonlin}
 \left\{
 \begin{array}{rcl}
F(\ep D^2 u^\ep,  D \ue, x,\xoe) +\ep  \ue & = & \e f(x, \xoe)   
\quad \hbox{in  }\Oe \, ,\\
\noalign{\vskip6pt}
 L ( D\ue ,x, \xoe  ) & = & 0 \quad \hbox{on  }\dOe\; ,
  \end{array}
 \right.
\end{equation}
and assume that 
$F$ and $L$  satisfy (F0), (F1), (F2), (F3), (L1) and (L3). 
Note that,
compared to the linear setting here, we have multiplied the equation by $\e$ 
here to simplify the presentation.

\medskip
Our assumption on the first cell problem is that 
$$\text{(H4)}\qquad \left\{ \enspace
\vcenter{\baselineskip=20pt\vsize=50pt\hsize=5.0truein \noindent
for all $x,p\in\real^N$, there exists a bounded, periodic in $y$ 
solution $\bar v(p,x,y)$ of
\newline 
\centerline{$F(D_{yy}^2 \bar v, p + D_y \bar v,x,y) =0\quad\hbox{in  }\R\; ,$}
\newline
which is 
a smooth function with respect to all its variable.
} \right.  \hskip.5truein $$

\medskip
Next we assume that $F$ is smooth in $M$ and $p$. 
Using the regularity of  $\bar v$ in (H4) we set 
$$ A(p,x,y)= F_M (D_{yy}^2\bar v,  p + D_y\bar v,   x,y ) 
\  \text{ and }\  
b(p,x,y) = F_p (D_{yy}^2\bar v,  p + D_y  \bar v,x,y )\; .$$

The cell problem for the second corrector is:  
\begin{equation*}
\text{(H5)}\quad \left\{ \enspace
\vcenter{\baselineskip=22pt \vsize=50pt\hsize=5.0truein \noindent
For all $x$, $p\in\real^N$ and $M\in S^N$, there exists
a unique constant $\Fb (M,p,x)$ such that the equation
\newline 
\centerline{$-\tr(A(x,p,y)\left(M+2 (D_{xy}^2 \bar v + MD_{yp}^2\bar v)
+ D_{yy}^2 \bar w \right)$ }
\newline
$- b(x,p,y) \cdot (D_{x}^2 \bar v + MD_{p}^2\bar v + D_y \bar w)
= f(x,y) - \Fb (M,p,x) \quad\hbox{in  } \R$  
\newline 
has a bounded, periodic in $y$ solution $\bar w$, 
which depends smoothly on $M$, $p$, $x$ and $y$.
} \right.  \hskip.9truein
\end{equation*}

\medskip
Finally, to specify the boundary condition, we assume that (E1) and 
(F6) are also satisfied, 
and we use the results of Section~\ref{sec:exist} and 
\ref{sec:uni}, which yield the existence of a unique constant 
$\bar \mu (p,x)$ such that the boundary ergodic problem
\begin{equation}\label{erghomnonlin}
 \left\{
 \begin{array}{rcl}
F(D_{yy}^2\tilde v,  p + D_y \tilde v,x,y) & = & 0 \quad\hbox{in  } \O \, ,\\
\noalign{\vskip6pt}
L (p+D_y \tilde v, x,y ) & = & -\bar \mu (p,x) \quad \hbox{on  }\dO\; .
  \end{array}
 \right.
\end{equation}
has a bounded, periodic in $y'$, solution $\tilde v$.  

The result is 

\begin{theor}\label{nonlinhom} 
Assume {\rm (F0), (F1), (F2), (F3), (F6), (E1),} 
{\rm (L1), (L3)},  {\rm (H4)}, {\rm (H5)} and that $F$ is a smooth 
function of $M$ and $p$. 
If $\bar \mu$ satisfies {\rm (L1)} and {\rm (L2)}, 
the sequence $(\ue)_{\e >0}$ converges, locally uniformly, 
as $\e\to0$, to the unique solution $\bar u$ of
\begin{equation*}	
 \left\{
 \begin{array}{rcl}
\Fb (D^2 \bar u,D\bar u,x) + \bar u = 0 \quad\hbox{ in  }\ \{x_N >0\} \, ,\\
\noalign{\vskip6pt}
\bar \mu (D\bar u,x) = 0 \quad \hbox{ on  }\ \{x_N = 0\}\; .
  \end{array}
 \right.
\end{equation*}
\end{theor}

We omit the proof since it follows very closely 
the proof of Theorem~\ref{th:homlin}. 
We remark that the nonlinear boundary condition does not introduce 
any additional difficulties. 
The most restrictive assumption concerns the smoothness of $F$ in $M$ and 
$p$, which is not easily satisfied by Hamilton-Jacobi-Bellman equations. 
At least Theorem~\ref{nonlinhom} applies to cases where $F$ 
is of the form $F= F_1 +H$, where $F_1$ a linear second order operator and $H$ 
is nonlinear in $p$ and smooth in $p$ and $x$.

\section{Open Problems and Remarks}

$\quad$ 
The periodicity of the equation and the boundary is a 
key assumption in our approach. 
The existence of the boundary cost and of the associated solution of the 
Neumann problem uses this property in a very strong way. 

It would be interesting to know whether approximate correctors exists  when 
the equation and the boundary condition are periodic 
but the  half-space domain does not have 
the right periodicity property, i.e., it is typically a hyperplane 
with an irrational slope. 
This seems to lead to a framework rather close to the one considered by 
Ishii \cite{I1}, but at the moment we are unable to obtain results 
in this direction.
Of course the ultimate goal is to consider the homogenization problem 
in random environments.

The following remark illustrates one of the difficulties to obtain 
such type of result. 
For $q\in \real^N\setminus \{0\}$, consider the half-space 
$$H_q=\{x \in \R\, : q \cdot x > 0\}$$ 
and assume that the problem 
\begin{equation*}
\left\{
\begin{array}{c}
-\Delta u = 0 \quad \hbox{in  }H_q\; ,\\
D u\cdot n 
= g(x) + \mu \quad \hbox{on  }\partial H_q\, ,
\end{array}\right.
\end{equation*}
has a bounded solution $u$ for $\hat a$ continuous, $\Z^N$-periodic. 
Then it is easy to see that we must have 
$$ \mu (q) = - \lim_{R \to + \infty}\, |B(0,R)\cap H_q|^{-1}\int_{B(0,R)\cap H_q}\, g(x) dx\, .$$

Now, in $\r^2$, choose $g(x_1,x_2)= \tilde g (x_2)$ where $\tilde g$ is 
a $\Z$-periodic function. 
It follows $\mu(e_2) = \tilde g (0)$, while, 
if $q=e_2+ \alpha e_1$ with $\alpha$ is a small parameter, 
$\mu(e_2) =  \int_0^1 \, \tilde g(s) ds$. 
This shows that a priori $\mu(q)$ is not a continuous function of $q$. 
Therefore an argument by approximation of a non-periodic situation by a 
periodic one cannot be used to prove the existence of $\mu$ and of 
the approximate corrector in the non-periodic framework.


\begin{thebibliography}{99}



\bibitem{A3} Arisawa, M.: 
{\sl~Long time averaged reflection force and homogenization of oscillating 
Neumann boundary conditions,}  
Ann. Inst. H. Poincar\'e Anal. Lin\'eaire 20 (2003)  293--332.


\bibitem{bb} Barles, G.: 
{\sc Solutions de viscosit\'e des \'equations de Hamilton-Jacobi}. 
Collection ``Math\'ematiques et Applications'' de
la SMAI, n$^\circ$17, Springer-Verlag (1994).

\bibitem{b2}
Barles, G.:  
{\sl Nonlinear Neumann boundary conditions for
quasilinear, degenerate elliptic equations and applications,}
Journal of Diff. Eqs. { 154} (1999) 191--224.

\bibitem{BDL-eb} Barles, G. and  Da Lio, F.: 
{ \sl On the boundary ergodic problem for fully nonlinear equations 
in bounded domains with general nonlinear Neumann boundary conditions}, 
Ann. IHP, Analyse non lin\'eaire, to appear.

\bibitem{BDL} 
Barles, G.   and   Da Lio, F.: 
{\sl Local $C^{0,\alpha}$ estimates for viscosity solutions of 
Neumann-type boundary value problems}, preprint.

\bibitem{BaLi}  Barles, G. and  Lions, P.-L.: 
{\sl Remarques sur les probl\`emes de reflexion obliques,} 
C. R. Acad. Sci. Paris, t. 320 Serie I (1995)  69--74.


\bibitem{BS1} Barles, G. and Souganidis, P.E.: 
{\sl On the large time behaviour of solutions of Hamilton-Jacobi equations,} 
SIAM J. Math.  Anal.  31 (2000) 925--939.

\bibitem{BS2} Barles, G. and Souganidis, P.E.: 
{\sl Space-time periodic solutions and long-time behavior of solutions 
to quasi-linear parabolic equations,}  
SIAM J. Math. Anal. 32 (2001) 1311--1323.

\bibitem{BLP} 
Bensoussan, A., Lions, J.-L., and Papanicolaou, G.: 
{\tt Cited on page 3.}

\bibitem{B} Bensoussan, A.: 
{\sc Perturbation methods in optimal control,}
Translated from the French by C. Tomson. Wiley/Gauthier-Villars Series 
in Modern Applied Mathematics.  
John Wiley \& Sons, Ltd., Chichester; Gauthier-Villars, Montrouge, 1988.


\bibitem{cil}
Crandall M.G., Ishii, H. and  Lions, P.-L.: 
{\sl User's guide to viscosity solutions of second order Partial differential
equations}, Bull. Amer. Soc. {27} (1992)  1--67.


\bibitem{E1} Evans, L.C.: 
{\sl The perturbed test function method for viscosity solutions of 
nonlinear PDE,}  
Proc. Roy. Soc. Edinburgh Sect. A 111 (1989)  359--375.

\bibitem{E2} Evans, L.C.: 
{\sl Periodic homogenisation of certain fully nonlinear partial differential 
equations,} 
Proc. Roy. Soc. Edinburgh Sect. A 120 (1992)  245--265.

\bibitem{I1} Ishii, H.: 
{\sl Almost periodic homogenization of Hamilton-Jacobi equations,} 
International Conference on Differential Equations, 
Vol. 1, 2 (Berlin, 1999), 600-605, World Sci. Publishing, River Edge, NJ, 2000.

\bibitem{I2} Ishii, H.: 
{\sl Perron's method for Hamilton-Jacobi Equations,} 
Duke  Math. J. {55} (1987) 369--384.




\bibitem{LS1} Lions P.-L and Souganidis, P.E.:
viscous

\bibitem{LS2} Lions P.-L and Souganidis, P.E.:
periodic

\bibitem{LiSz} Lions P.-L and Sznitman, A.S.: 
{\sl Stochastic differential equations with
reflecting boundary conditions, }
Comm. Pure and Applied Math.  37	
(1984) 511--537. 

\bibitem{HT} Tanaka, H.: 
Homogenization of diffusion processes with boundary conditions. 
Stochastic analysis and applications, 411--437, Adv. Probab. 
Related Topics 7 Dekker, New York, 1984.

\end{thebibliography}
\end{document}